\documentclass[a4paper,twocolumn]{article} 

\usepackage[utf8]{inputenc}
\usepackage{amssymb,amsmath,amsthm} 
\usepackage{caption, subcaption} 
\usepackage{epsfig} 
\usepackage[numbers]{natbib} 
\usepackage{url} 

\usepackage{graphicx}      
\usepackage{psfrag}
\usepackage{algorithm}
\usepackage{import}
\usepackage{algpseudocode}
\usepackage{color}
\usepackage{booktabs}

\newcommand{\R}{\mathbb{R}}

\title{Accelerated Nonlinear Model Predictive Control by Exploiting Saturation} 

\author{Raphael Dyrska, Ruth Mitze, and Martin M\"onnigmann\\
	Automatic Control and Systems Theory, Dep. of Mechanical Engineering,\\
	Ruhr-Universit\"at Bochum, 44801 Bochum, Germany.\\ E-mail: {\tt\small raphael.dyrska@rub.de}, {\tt\small ruth.mitze@rub.de}, and {\tt\small martin.moennigmann@rub.de}}

\begin{document}
\maketitle

\begin{abstract}                
We present an approach for accelerating nonlinear model predictive control. 
If the current optimal input signal is saturated, also the optimal signals in subsequent time steps often are.
We propose to use the open-loop optimal input signals whenever the first and some subsequent input signals are saturated. We only solve the next optimal control problem, when a non-saturated signal is encountered, or the end of the horizon is reached.
In this way, we can save a significant number of NLPs to be solved while on the other hand keep the performance loss small. Furthermore, the NMPC is reactivated in time when it comes to controlling the system safely to its reference.

\end{abstract}


\section{Introduction}
Model predictive control (MPC) is based on recurrently solving an optimal control problem (OCP) and using the first resulting optimal control input for the current system state. 
Research on MPC has focused on reducing the time required to solve the OCP or even to avoid its solution at all.
The first topic is either addressed by algorithmic approaches~\citep{Best1996, Ferreau2008, Diehl2009}, or for the linear case, by approaches that exploit the structure of the OCP~\citep{Jost2017}.

In contrast, there exist (again mostly for the linear case) a lot of approaches to avoiding solution of an OCP. These approaches are based on the explicit solution~\citep{Bemporad2002, Tondel2003} and variants thereof~\citep{Johansen2004, Kvasnica2015}, or on the insights regarding the piecewise affine control laws provided by the explicit solution~\citep{Jost2015a, BernerP2018}.

Nonlinear MPC (NMPC) is more challenging, since a lot of assumptions and helpful properties such as convexity do no longer hold, and solving a nonlinear program (NLP) may be very expensive~\citep{Azhmyakov2008, Lautenschlager2015}. There exist some extensions of explicit methods for linear systems to the nonlinear case~\citep{Johansen2002, SchulzeDarup2012_NMPC, Monnigmann2015a, Oberdieck2016}, but none of them are as comprehensive as in the linear case. 

Large areas of the feasible state space, i.e., the set of initial states for which a solution satisfying all constraints exists, result in saturated control inputs. 
In this paper, we propose a method that avoids the solution of some NLPs based on information about saturated control inputs. The proposed method is not based on any explicit solution.

It is well known that for deterministic linear MPC, assuming that the terminal set and the terminal cost are chosen appropriately, the predicted open-loop optimal solution and the resulting closed-loop optimal solution are equal, if the terminal constraints are inactive~\citep{Chmielewski1996, Monnigmann2019}.
This statement does no longer hold for NMPC, among other reasons because the unconstrained optimal feedback law and the terminal set can typically only be approximated in the nonlinear case~\citep{Chen1998, Mayne2000}. 
Nevertheless, it can be observed that the open-loop optimal input signals are often closed-loop optimal if the first and some subsequent signals are saturated.

We present a simple but effective heuristic method of avoiding the solution of NLPs: If the OCP for the current state results in saturated control inputs, we reuse this input as long as it is predicted to be open-loop optimal for future time steps. 
In the example section, we will see that this simple heuristic results in skipping some NLPs, while leading to almost the same results as a standard NMPC scheme.

\section{Problem statement and\\ notation}\label{sec:ProblemStatement}

Throughout the paper, we consider discrete-time, nonlinear systems of the form
\begin{align}\label{eq:system}
  x(k+1)=f(x(k),u(k)),\, k= 0, 1, \dots
\end{align}
with $f\in\mathbb{C}^2$, $f(0,0)=0$.
States and inputs are assumed to be constrained according to
\begin{align*}
  x(k) &\in \mathcal{X}\subset\R^n,
  \\
  u(k) &\in \mathcal{U}\subset\R^m
\end{align*}
for all time steps $k$, where 
\begin{align}
\mathcal{X}&=\{x: \underline{x}\leq x\leq\overline{x}\},\\
\mathcal{U}&=\{u: \underline{u}\leq u\leq\overline{u}\},
\end{align}
for some $\underline{x}, \overline{x}\in\R^{n}$, $\underline{u}, \overline{u}\in\R^{m}$ such that $0\in \text{int}(\mathcal{X})$ and $0\in \text{int}(\mathcal{U})$. 

To regulate system~\eqref{eq:system} to the origin, we solve, on a receding horizon, the optimal control problem
\begin{subequations}\label{eq:OCP}
\begin{align}\label{eq:OCPCost}
 V_{N}(x)=\ \underset{X, U}{\text{min}}\ V_f(x(N))+ \sum_{k=0}^{N-1} \ell(x(k), u(k))
\end{align}
subject to 
\begin{align}\label{eq:OCPConstraints}
  \begin{split}
    &x(k+1)=f(x(k), u(k)),\ k=0,...,N-1\\
    &x(k) \in \mathcal{X},\ k=0,...,N-1\\
    &u(k) \in \mathcal{U},\ k=0,...,N-1\\
    &x(N) \in \mathcal{T},
  \end{split}
\end{align}
\end{subequations}
for given initial value $x(0)$, states $X=(x^{T}(1),...,x^{T}(N))^{T}$, inputs $U=(u^{T}(0),...,u^{T}(N-1))^{T}$, prediction horizon $N$, terminal cost
$V_f(x)= x^T P x$, and stage cost $\ell(x, u)= x^T Q x+ u^T Ru$. The weighting matrices $P$, $Q$, and $R$ are assumed to be positive definite matrices of the obvious dimensions. The terminal set $\mathcal{T}\subseteq\mathcal{X}$ with $0\in\text{int}(\mathcal{T})$ and the terminal cost function fulfill, together with a stabilizing control law $\kappa(x)$ on $\mathcal{T}$, the requirements of a stabilizing triple as defined in~\cite{SchulzeDarup2012_NMPC} (see also~\cite{Chen1998, Mayne2000} for stability properties).

For simplicity, we rewrite the optimal control problem~\eqref{eq:OCP} by substituting the system dynamics \eqref{eq:system} into \eqref{eq:OCP}
\begin{subequations}\label{eq:OCPCompact}
  \begin{align}\label{eq:OCPCost_compact}
    &\underset{U}{\text{min}}\ V(x(0), U)
  \end{align}
  subject to
  \begin{align}\label{eq:constraintsCompact}
    G(x(0), U)\leqslant 0.
  \end{align}
\end{subequations}
Note that both problems describe the same NLP. 

\subsection*{Notation}
Let the set of all states for which problem~\eqref{eq:OCPCompact} (or equivalently~\eqref{eq:OCP}) has a solution be denoted $\mathcal{F}$. For any $x(0)\in\mathcal{F}$, let $V(x(0),U^{*}(x(0)))$ refer to the optimal solution of~\eqref{eq:OCPCompact} with optimizer $U^{*}(x(0))$. A constraint $i$ is called active for $x(0)\in\mathcal{F}$ if
$G_{i}(x(0), U^{*}(x(0)))=0$, where the index $i$ refers to row $i$ of $G$. Let $\mathcal{A}(x(0))$ refer to the set of active constraints for $x(0)$. 
If $u(0)$ is subject to the bounds $u(0)\in [\underline{u}_{0}, \overline{u}_{0}]$, we call $u(0)$ saturated whenever $u(0)=\underline{u}_{0}$ or $u(0)=\overline{u}_{0}$. 
We say e.g. '$\overline{u}_{k}$ is active' for any constraint and any time step $k$ if the index $i$ referring to the corresponding line in matrix~\eqref{eq:constraintsCompact} is in the active set, i.e., $i\in\mathcal{A}$.

\section{Saturated control laws}\label{sec:method}

In model predictive control, problem~\eqref{eq:OCP} is recurrently solved on a receding horizon, and in every time step, the first optimal input $u^{*}(0)$ is applied to the system. 

\begin{figure}[H]
\begin{center}
\input{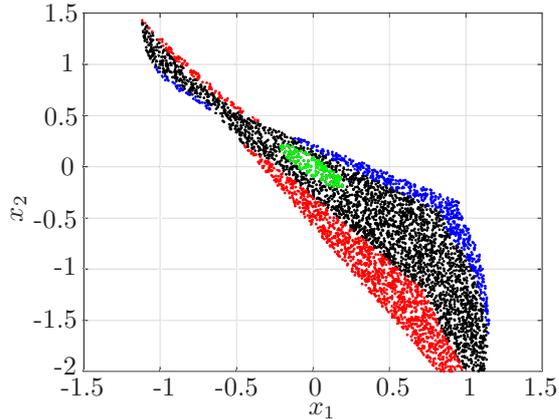}
\caption{Feasible set of the system~\eqref{eq:example}. For red (blue, black) states the upper (lower, no) bound on $u(0)$ is active. Green states are part of the terminal set.}
\label{fig:feasibleSet}
\end{center}
\end{figure}

The optimal input frequently is saturated close to the boundary of the feasible set.
For optimization problems with tight input constraints $\mathcal{U}$, the saturated regions can easily make up $40\%$ of the whole feasible set. Figure~\ref{fig:feasibleSet} shows this for the example from Section~\ref{sec:example} (see also Table~\ref{tab:distribution}). 
For initial values lying in these regions, some of the following predicted open-loop, as well as the actual used closed-loop control inputs, also take on saturated values. This effect is larger for the open-loop solution since, roughly speaking, the controller needs to reach the terminal set in $N$ time steps.
Often, however, the open-loop and closed-loop optimal inputs attain the same saturated values for some time steps.

\begin{table}[b]
\begin{center}
\caption{Distribution of initial values resulting in inputs at the lower ($\underline{u}_{0}$) and upper ($\overline{u}_{0}$) bound and other inputs in percent.}
\begin{tabular}{cccc}\label{tab:distribution}
resulting control input & $\underline{u}_{0}$ & $\overline{u}_{0}$ & other \\\hline
percentage of $\mathcal{F}$ & $12.5$ & $28.18$ & $59.32$ \\ \hline
\end{tabular}
\end{center}
\end{table}

\subsection{Reusing predicted control laws}
It is our central idea to reuse the saturated control law as long as it was predicted for the initial value. When solving the NLP for an initial value $x_{0}$, we check (i) if the resulting input $u(0)$ is saturated and (ii) if there are subsequent time steps for which the corresponding constraints on the input are predicted to be active as well. If both conditions hold, instead of solving the NLP, we reuse the saturated control law as long as it was predicted. After all predicted inputs were used, we switch back to the NMPC in the classic fashion.

The proposed procedure is summarized in Algorithm~\ref{alg:Algorithm}. As mentioned before, after solving the problem~\eqref{eq:OCPCompact}, the algorithm basically checks the constraints regarding $u(0)$ and, in case one of them is active, also checks the predicted input values for $\tilde{k}>0$. If at least one prediction is saturated as well, instead of solving the NLP in the subsequent time steps, the saturated values are used as long as they were predicted for the initial state. Whenever the end of this prediction is reached, classic NMPC is used to steer the state closed-loop to the origin. Note here that this procedure is only performed at the beginning of the NMPC, i.e., only when solving for the first time step.

While the procedure described in Algorithm~\ref{alg:Algorithm} only applies to systems with a single input, i.e., $m=1$, it can be easily extended to the multiple input case. However, all inputs have to be saturated in this case. The open-loop solution can then be used until the shorter number of saturated time steps is reached.

\begin{algorithm}[t]
\begin{algorithmic}[1]
\caption{NMPC with reusing predicted saurated control inputs}\label{alg:Algorithm}
\State \textbf{Input:} $\mathcal{A}(x(0))$.
	\If{$\overline{u}_{0}$ is active}
		\If{$\overline{u}_{\tilde{k}}$ is active for $\tilde{k}\geq 1$}
		\State $u(1,...,\tilde{k})=u(0)=\overline{u}$.
		\EndIf
	\ElsIf{$\underline{u}_{0}$ is active}
		\If{$\underline{u}_{\tilde{k}}$ is active for $\tilde{k}\geq 1$}
		\State $u(1,...,\tilde{k})=u(0)=\underline{u}$.
		\EndIf
	\Else
	\State solve~\eqref{eq:OCPCompact} for $k\geq 1$.
	\EndIf
\end{algorithmic}
\end{algorithm}

\section{Example}\label{sec:example}

We compare the proposed approach to classic NMPC in terms of the number of solved NLPs. We analyze both the nominal and additively disturbed case. 

\subsection{Simulation setup}

We consider the nonlinear discrete-time benchmark example from~\cite{SchulzeDarup2012_NMPC}
\begin{small}
\begin{align}\label{eq:example}
\begin{split}
x_{1}(k+1)&=x_{1}(k)+0.1x_{2}(k)+0.1(0.5+0.5x_{1}(k))u(k),\\
x_{2}(k+1)&=x_{2}(k)+0.1x_{1}(k)+0.1(0.5-2.0x_{2}(k))u(k),
\end{split}
\end{align}
\end{small}
with state and input constraints
\begin{align*}
  -2\leqslant x_{1, 2}(k)\leqslant 2,
  \\
  -0.5\leqslant u(k)\leqslant 0.5.
\end{align*}
The weighting matrices are chosen as
\begin{align*}
Q=\begin{pmatrix}
0.05 & 0\\
0 & 0.05
\end{pmatrix},\ \ \ \ R=0.1.
\end{align*}
The terminal penalty matrix
\begin{align*}
P=\begin{pmatrix}
5.9353 & 5.2774\\
5.2774 & 5.9353
\end{pmatrix}
\end{align*}
and the terminal set $\mathcal{T}=\lbrace x\in \mathcal{X}| x^{T}Px\leq \alpha\rbrace$, $\alpha=0.0606$ 
fulfill the requirements to form a stabilizing triple. We used a prediction horizon of $N=12$.

We simulated the regulation to the terminal set in the classic fashion and using the proposed method 
for the $5000$ random feasible initial values shown in Fig.~\ref{fig:feasibleSet}. 
Almost $12\%$ of the 5000 initial values resulted in an optimal control input $u(0)$ saturated at the lower bound. For the upper bound saturation resulted for about $28\%$ of all cases. We counted the number of solved NLPs and measured the cost performance $\tilde{V}$ in terms of
\begin{align*}
\hat{V}=\sum_{k=0}^{\hat{k}}x(k)^{T}Qx(k)+u(k)^{T}Ru(k),
\end{align*}
where $\hat{k}$ is the number of time steps needed to reach the terminal set. 
In other words, we calculated the weighted distance of the current state and input from the origin in every time step. 

\begin{figure}[t]
\input{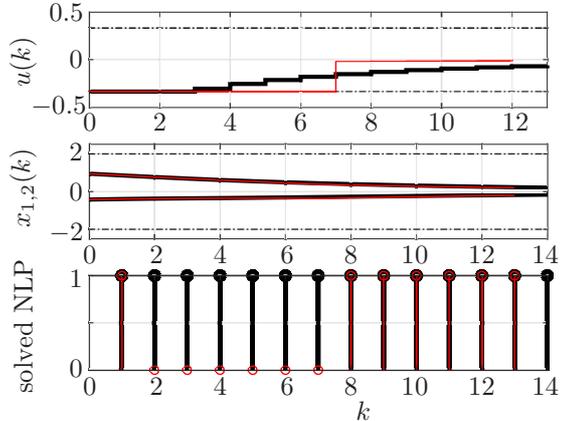}
\caption{Input and state trajectories and an indicator for solved NLPs for standard NMPC (black) and NMPC using the proposed method (red) for one of the 5000 initial values over time steps $k$. 
A value of $1$ (respectively $0$) indicates an (respectively no) NLP was solved (bottom plot).}
\label{fig:example}
\end{figure}

\begin{table}[b]
\begin{center}
\caption{Results for lower bound $\underline{u}$}
\begin{tabular}{cccc}\label{tab:lower}
 & classic & here & percent \\\hline
cost\\ Performance & $0.377$ & $0.386$ & $2.3873$ \\
no. of\\ solved NLPs & $12.8688$ & $10.6752$ & $17.033$ \\ \hline
\end{tabular}
\end{center}
\end{table}

\begin{table}[b]
\begin{center}
\caption{Results for upper bound $\overline{u}$}
\begin{tabular}{cccc}\label{tab:upper}
 & classic & here & percent \\\hline
cost\\ Performance & $0.4407$ & $0.442$ & $0.295$ \\
no. of\\ solved NLPs & $12.1001$ & $8.5578$ & $29.275$ \\ \hline
\end{tabular}
\end{center}
\end{table}

\subsection{Deterministic case}

Figure~\ref{fig:example} shows the resulting state and input trajectories as well as the number of solved NLPs for one of the initial values for the deterministic case. In this illustrative example, problem~\eqref{eq:OCPCompact} resulted in a $u(0)$ saturated at the lower bound, $u(0)=-0.5$. As $u=-0.5$ was predicted for the next six time steps as well, we reuse the open-loop solution and therefore do not have to solve the corresponding NLPs during these time steps using the proposed method. Note that we even reach the terminal set one time step earlier compared to the classic approach because the open-loop solution is naturally more aggressive.

For all among the $5000$ initial values resulting in an input at lower (upper) bound, Table~\ref{tab:lower} (Table~\ref{tab:upper}) compares the mean values of the cost performance and the number of solved NLPs resulting for classic and modified NMPC.

For initial values resulting in a saturated input at the lower bound (blue states in Fig.~\ref{fig:feasibleSet}) we save on the average around $17\%$ of NLPs. The cost performance, however, drops only by around $2.3\%$.

For initial values resulting in an input saturated at the upper bound, we save almost $30\%$ by using the proposed method, while losing a negligible amount of $0.3\%$ regarding the cost performance.

\begin{figure}[t]
\input{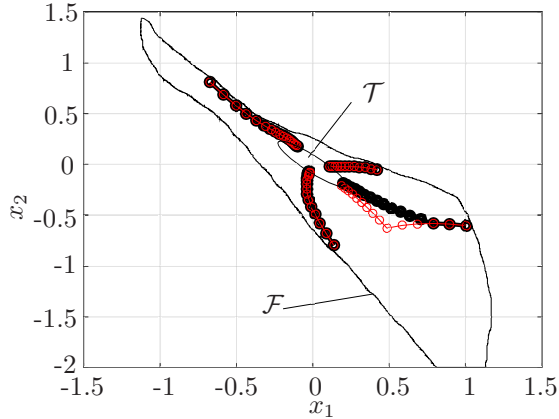}
\caption{Trajectories of NMPC (black) and NMPC reusing predicted saturated control inputs (red) over the approximated feasible set $\mathcal{F}$ and terminal set $\mathcal{T}$.}
\label{fig:feasibleSetTrajectories}
\end{figure}

In conclusion, concerning the whole feasible set (i.e., for all initial values, regardless if saturated input or not), we can save on the average around $10.38\%$ of NLPs, while the control performance drops only very midly. Figure~\ref{fig:feasibleSetTrajectories} shows illustrative state trajectories resulting from classic NMPC (black) and the proposed method (red), respectively, on an approximation of the feasible set presented in Figure~\ref{fig:feasibleSet}. While all trajectories are very similar, the trajectory starting at $x_{0}=(1.004, -0.6015)$ shows that the number of predicted saturated control inputs is higher than the number used in the closed-loop fashion. 

\subsection{Disturbed case}

We also evaluated the method under the presence of bounded disturbances. The dynamic model used to determine the subsequent state  reads
\begin{align}
x(k+1)=f(x(k), u(k))+w,
\end{align}
with disturbance vector $w\in\mathcal{W}$. Constraints $\mathcal{W}$ are chosen such that
\begin{align*}
-0.01\leq w\leq 0.01.
\end{align*}
We used the same initial values as for the deterministic case for which a saturated control input minimized~\eqref{eq:OCPCompact} and applied the same random disturbance vector to the trajectory resulting from each initial value with and without the proposed method. No robustification method was applied, and $7$ out of the $2034$ initial values resulted in infeasible behavior. They were removed from the statistics.

For the remaining initial values, Table~\ref{tab:lower_dist} (Table~\ref{tab:upper_dist}) compares again the number of NLPs and the cost performance for the lower (upper) bound under the presence of disturbances. About $16.6\%$ of the NLPs were skipped for initial values resulting in a saturated input at the lower bound. About $29\%$ were skipped for initial values with inputs saturated at the upper bound. The cost performances are again negligibly small. 

\begin{table}[H]
\begin{center}
\caption{Results for lower bound $\underline{u}$ and present disturbances}
\begin{tabular}{cccc}\label{tab:lower_dist}
 & classic & here & percent \\\hline
cost\\ Performance & $0.3799$ & $0.3892$ & $2.3895$ \\
no. of\\ solved NLPs & $13.0516$ & $10.8806$ & $16.634$ \\ \hline
\end{tabular}
\end{center}
\end{table}

\begin{table}[H]
\begin{center}
\caption{Results for upper bound $\overline{u}$ and present disturbances}
\begin{tabular}{cccc}\label{tab:upper_dist}
 & classic & here & percent \\\hline
cost\\ Performance & $0.4425$ & $0.4438$ & $0.2929$ \\
no. of\\ solved NLPs & $12.2011$ & $8.6468$ & $29.131$ \\ \hline
\end{tabular}
\end{center}
\end{table}

\section{Conclusion}\label{sec:conclusion}
We presented a simple but effective method for accelerating nonlinear model predictive control by avoiding the solution of potentially expensive NLPs. 
Saturated regions can make up a large part of the feasible set. 
We were able to show for both, the deterministic and non-deterministic case, that using a simple heuristic results in a reduction of NLPs, with a negligible loss of control performance. 
While the treated example was very simple and had only one input, the result warrant further investigations in the proposed method. 

\section{Acknowledgements}
Support by the Deutsche Forschungsgemeinschaft (DFG) under grant MO 1086/15-1 is gratefully acknowledged.

\bibliographystyle{plain}
\bibliography{msdLiterature}

\end{document}